\documentclass{lmcs}
\pdfoutput=1

\usepackage[utf8]{inputenc}

\usepackage{lastpage}
\lmcsdoi{18}{4}{13}
\lmcsheading{}{\pageref{LastPage}}{}{}%
{Mar.~29,~2021}{Dec.~20,~2022}{}

\usepackage{enumitem}
\setenumerate{label=(\roman*)}

\usepackage{thmtools}

\usepackage{graphicx}
\usepackage{amsfonts}
\usepackage{amsmath}
\usepackage{amssymb}

\usepackage{tikz-cd}

\usepackage{bbm}
\usepackage{stmaryrd} 

\usepackage[colorlinks=false,hidelinks]{hyperref}

\newcommand{\heyting}[1]{{\llbracket #1 \rrbracket}}
\newcommand{\set}[1]{\left\{ #1 \right\}}
\newcommand{\Set}[2]{\left\{ #1 \, : \, #2 \right\}}
\renewcommand{\phi}{\varphi}
\renewcommand{\epsilon}{\varepsilon}
\renewcommand{\P}{\mathcal{P}}
\newcommand{\apart}{\mathrel{\#}}

\newcommand{\numeral}[1]{\bar{ #1 }}

\newcommand{\inl}{\mathsf{inl}}
\newcommand{\inr}{\mathsf{inr}}
\newcommand{\fst}{\mathsf{fst}}
\newcommand{\snd}{\mathsf{snd}}
\newcommand{\Nsucc}{\mathsf{succ}}
\newcommand{\pair}{\mathsf{pair}}
\DeclareMathOperator{\Pow}{Pow}
\DeclareMathOperator{\Sub}{Sub}

\newcommand{\Asm}{\mathbf{Asm}}
\newcommand{\ApAsm}{\mathbf{ApAsm}}
\newcommand{\ApAss}{\ApAsm}
\newcommand{\ApType}{\mathbf{ApType}}
\renewcommand{\Im}{\mathrm{Im}}

\newcommand{\CE}{\mathsf{CE}}
\newcommand{\AC}{\mathsf{AC}}
\newcommand{\IP}{\mathsf{IP}}
\newcommand{\MP}{\mathsf{MP}}
\newcommand{\EXT}{\mathsf{EXT}}

\newcommand{\usftext}[1]{\textsf{\upshape #1}}
\newcommand{\eha}{\ensuremath{{\usftext{E-HA}}^\omega}} 

\newcommand{\tca}[1]{\mathcal{#1}}

\newcommand{\IL}{\ensuremath{\usftext{IL}}}


\begin{document}

\title[Converse extensionality and apartness]{Converse extensionality and apartness}

\author{Benno van den Berg\lmcsorcid{0000-0002-0469-0788}}
\author{Robert Passmann\lmcsorcid{0000-0002-7170-3286}}
\address{Institute for Logic, Language and Computation, Universiteit van Amsterdam, Postbus 90242, 1090 GE Amsterdam, The Netherlands}
\email{B.vandenBerg3@uva.nl, r.passmann@uva.nl}

\thanks{Robert Passmann was supported by a doctoral scholarship of the \emph{Studien\-stiftung des deutschen Volkes}.}


\begin{abstract}
    
\noindent In this paper we try to find a computational interpretation for a strong form of extensionality, which we call ``converse extensionality''. Converse extensionality principles, which arise as the Dialectica 
interpretation of the axiom of extensionality, were first studied by Howard. In order to give a computational interpretation to these principles, we reconsider Brouwer's apartness relation, a strong constructive form of inequality. Formally, we provide a categorical construction to endow every typed combinatory algebra with an apartness relation. We then exploit that functions reflect apartness, in addition to preserving equality, to prove that the resulting categories of assemblies model a converse extensionality principle.
\end{abstract}

\maketitle

\section{Introduction}

Following Kreisel one of the main concerns of proof theory has become the extraction of hidden computational information from proofs. For this purpose G\"odel's Dialectica interpretation (combined with negative translation, if necessary) has proven itself to be indispensable. Indeed, within proof mining functional interpretations of various kinds have become a sophisticated and flexible tool for extracting additional qualitative and quantitative information from proofs (see \cite{kohlenbach08}).

One of the hardest principles to interpret using a functional interpretation is the principle of function extensionality. This principle, which says that two functions are equal if they yield the same output on the same input, is pervasive in mathematics. But it has proven difficult to interpret using the Dialectica interpretation, the reason being that the Dialectica interpretation requires one to interpret a stronger form of extensionality, which we have dubbed \emph{converse extensionality}:
\[ \CE_n: \quad \exists X \, \forall \Phi^{n+2} \, \forall f, g \, \big( \, \Phi f \not=_0 \Phi g \to f (X \Phi f g) \neq_0 g (X \Phi f g) \, \big). \]
Note that this is equivalent to
\[ \exists X \, \forall \Phi^{n+2} \, \forall f, g \, \big( \,  f (X \Phi f g) =_0 g (X \Phi f g) \to \Phi f =_0 \Phi g \, \big) \]
since equality of type 0 is decidable. As shown by Howard (see \cite[Appendix]{Troelstra344}), $\CE_0$ cannot be witnessed in the term model of G\"odel's $T$ and $\CE_1$ is unprovable in Zermelo-Fraenkel set theory (without choice). This has often been taken as an indication that a computational interpretation of function extensionality is well-nigh impossible.

The starting point for this paper was the question whether the situation is really that hopeless. Our idea is that by a suitable enrichment of data it might still be possible to interpret (fragments of) converse extensionality. For this we are looking at Brouwer's notion of \emph{apartness}.

Brouwer's idea was that equality might not be a primitive concept and could be defined as the negation of a strong notion of inequality called apartness. The paradigmatic example are the real numbers, where two reals $r$ and $s$ are apart when there are disjoint intervals with rational endpoints $I_1$ and $I_2$ such that $r \in I_1$ and $s \in I_2$. Equality of real numbers can then be defined as not being apart. The notion of apartness has continued to play a role in constructive mathematics to this very day (see \cite{darpomitrovic21} for a recent example).

Typical properties of the apartness relation $\apart$ are the following:
\begin{displaymath}
\begin{array}{c}
 \lnot x \apart x \\
  x \apart y \Rightarrow y \apart x \\
  x \apart y \Rightarrow (z \apart x \lor z \apart y)
\end{array}
\end{displaymath}
We will refer to these properties as \emph{reflexivity}, \emph{symmetry} and \emph{transitivity}, because these axioms ensure that equality $x = y := \lnot x \apart y$ has said properties. 

Our first step is the observation (see also \cite{TroelstraVanDalen88ii}) that on all the finite types equality can indeed be defined as the negation of a suitable notion of apartness. But that means that one may require functionals $f$ of type $\sigma \to \tau$ to come equipped with additional data that explains how from evidence that $fx$ and $fy$ are apart one obtains evidence that $x$ and $y$ are apart. Our initial results do indeed suggest that by enriching functionals with this data one may interpret certain forms of converse extensionality, although the results are not (yet) as strong as we had hoped.

To formulate the results that we have obtained so far, we use the notion of a typed combinatory algebra (tca), basically a model of G\"odel's $T$. We show that from every tca (including the term model of G\"odel's $T$) one can define a new tca, which we have dubbed the \emph{apartness types}. Our main result is that by using modified realizability over these apartness types one can interpret $\CE_0$. This shows (\emph{pace} Howard) that it might still be possible to interpret $\CE_0$ using terms from G\"odel's $T$. To interpret stronger principles ($\CE_1$ and higher) we currently have to use tcas which satisfy suitable continuity principles. 

For proving our results, we have decided to formulate them in a categorical framework, using categories of assemblies. So in Section 2 of the paper we will recall the definition of the category of assemblies over a tca, following Longley \cite{longley99}. We will also discuss the internal logic of the assemblies over a tca there: as far as we are aware, these results have not appeared earlier in the literature, but will not surprise the experts. In Section 3 we will show that the assemblies over what we will call an extensional tca satisfy principles reminiscent of Kreisel's modified realizability. These results were first obtained by Mees de Vries \cite{devries17} for one specific extensional tca: here we show that they hold more generally for any extensional tca. Finally, in Section 4 we introduce the apartness types and show that in the category of assemblies over the apartness types the converse extensionality principle $\CE_0$ holds. We have also included an appendix which explains our results from a proof-theoretic perspective.

Finally, we note that all our results are constructively valid, unless explicitly noted otherwise.

\section{The logic of assemblies}

The purpose of this section is to recall the definition of a typed combinatory algebra (tca) and the category of assemblies over a typed combinatory algebra. These definitions are due to Longley \cite{longley99} and can also be found in \cite{LietzStreicher2002} and \cite{longleynormann15}. We deviate from these sources by making two small changes: first of all, we will only consider total combinatory algebras, because all the examples that we will be interested in in this paper are total. Secondly, we include in the type structure of a typed combinatory algebra both finite sum types and a unit type. Our reason for doing so is that this allows us to prove that the assemblies over a tca form a Heyting category, which will be the main result of this section.

\subsection{Typed Combinatory Algebras}

We will start by defining typed combinatory algebras, the total variant of the typed partial combinatory algebras as in \cite{longley99}.

\begin{defi}
    \label{Definition: tca}
    A \emph{typed combinatory algebra} (tca) consists of a set $\tca{T}$ of types with the following data:
    \begin{enumerate}
        \item binary operations $\times$, $\rightarrow$, $+$ on $\tca{T}$, and distinguished types $\bot, \top, N$, 
        \item a set $|T|$ of realizers for every $T \in \tca{T}$,
        \item a total application function $\cdot_{S,T}: |S \rightarrow T| \times |S| \to |T|$,
    \end{enumerate}
    such that for all $U, S, T \in \tca{T}$ there are elements 
    \begin{align*}
        & \mathsf{exf} \in |\bot \to S|, \quad \mathsf{t} \in |\top|, \\
        & \mathsf{k}_{S,T} \in |S \to T \to S|, \quad
        \mathsf{s}_{S,T,U} \in |(S \to T \to U) \to (S \to T) \to (S \to U)|, \\
        & \mathsf{pair}_{S,T} \in |S \to T \to S \times T|, \quad
        \mathsf{fst}_{S,T} \in |S \times T \to S|, \quad
        \mathsf{snd}_{S,T} \in |S \times T \to T|, \\
        & \mathsf{inl}_{S,T} \in |S \to S + T|, \quad
        \mathsf{inr}_{S,T} \in |T \to S + T|, \\
        & \mathsf{case}_{S,T,U} \in |(S \to U) \to (T \to U) \to (S + T \to U)|, \\
        & \mathsf{0} \in |N|, \quad
        \mathsf{succ} \in |N \to N|, \quad
        \mathsf{R}_S \in |S \to (N \to (S \to S)) \to (N \to S))|,
    \end{align*}
    satisfying the following conditions
    \begin{align*}
        &
        \mathsf{k} \, a \, b = a, \quad
        \mathsf{s} \, a \, b \, c = a \, c \, (b \, c), \\
        &
        \mathsf{fst} \, (\mathsf{pair} \,a \, b) = a, \quad
        \mathsf{snd} \, (\mathsf{pair} \, a \, b) = b, \\
        &
        \mathsf{case} \, a \, b \, (\mathsf{inl} \, x) = a \, x, \quad
        \mathsf{case} \, a \, b \, (\mathsf{inr} \, x) = b \, x, \\
        &
        \mathsf{R} \, a \, b \, \mathsf{0} = a, \quad
        \mathsf{R} \, a \, b \, (\mathsf{succ} \, n) = b \, n \, (\mathsf{R} \, a \, b \, n),
    \end{align*}
    for $a$, $b$, $c$, $x$ and $n$ of the corresponding types.
\end{defi}

\begin{rem} As is customary in the theory of combinatory algebras, we usually omit the application $\cdot$ and write $a \, b$ or $a(b)$ instead of $a \cdot b$. In fact, we already started to do so when formulating the equations that the combinators should satisfy in the previous definition. In addition, the convention is that $\cdot$ associates to the left, so that $a \, b \, c \, d$ has to be read as $((a \cdot b) \cdot c) \cdot d$. Moreover, in any tca a form of lambda abstraction is available, in a manner similar to ordinary combinatory algebras (for which, see \cite{vanOosten2008}).
\end{rem}

\begin{defi} \label{Definition:finite types} If $\tca{T}$ is a tca, we will refer to the smallest set of types in $\tca{T}$ containing $\top$ and $N$ and closed under $\to$ and $\times$ as the \emph{finite types} in $\tca{T}$. In addition, we will also use natural numbers to refer to specific finite types, with $0 := N$ and $n+1 := n \to N$.
\end{defi}

\begin{defi} \label{Definition:properties of tcas}
    A typed combinatory algebra will be called
    \begin{enumerate}
        \item \emph{consistent} if $\mathsf{0} \not= \mathsf{succ \, 0}$.
        \item \emph{standard} if the mapping $\mathbb{N} \to |N|$ obtaining by sending $n$ to  the numeral $\overline{n} = \mathsf{succ}^n \mathsf{0}$ is a bijection.
        \item \emph{extensional} if for all types $S$ and $T$ the mappings
        \begin{align*}
            |S \times T| \to |S| \times |T| & : x \mapsto (\mathsf{fst} \, x, \mathsf{snd} \, x)  \\
            |T \to S| \to |S|^{|T|} & : x \mapsto \lambda y. x \cdot y
        \end{align*}
        are injective, and $\mathsf{t}$ is the sole element of $|\top|$.
    \end{enumerate}
\end{defi}

\begin{rem}
    Note that the map $p: |S| \times |T| \to |S \times T|: (x,y) \mapsto \mathsf{pair} \, x \, y$ is always injective, as it is a section of $x \mapsto (\mathsf{fst} \, x, \mathsf{snd} \, x)$. For this reason we will often use $(x,y)$ as an abbrevation for $\mathsf{pair} \, x \, y$. Hence, if $\tca{T}$ is extensional, the map $p$ will actually be a bijection (but this will not be the case in general).
\end{rem}

As in the case of partial combinatory algebras, every tca admits some recursion theory, see, e.g., van Oosten's book \cite[Chapter 1]{vanOosten2008}. In particular, for every fixed $n \in \mathbb{N}$, we can code finite sequences $(x_0,x_1,\dots,x_n)$ of length $n$, and $i$-th projections $\mathsf{proj}_i$ using just $\pair$, $\fst$ and $\snd$.

\begin{exa} \label{Example: tcas}
    Examples of typed combinatory algebras are abundant.
    \begin{enumerate}
        \item Every partial combinatory algebra $\mathcal{A}$ gives rise to a tca by taking the powerset of $\mathcal{A}$ as the set of types, with $|X| = X$, and the operations appropriately defined. We can also restrict to those subsets of $\mathcal{A}$ that are inhabited: this also gives one a tca.  
        \item We consider Kleene's first algebra $\mathcal{K}_1$ as a tca in the sense of (i) of the previous example.
        \item Similarly, $\mathcal{K}_2^\mathrm{rec}$ is the tca obtained as the recursive submodel of Kleene's second algebra $\mathcal{K}_2$.
        \item The closed terms of Gödel's $T$ form a tca, provided we take a version of G\"odel's $T$ which includes finite sum and unit types. The types are the types of G\"odel's $T$ and the realizers of a type consist of the closed terms of that type. This shows that unbounded search is generally not available in tcas.
        \item If $\mathcal{C}$ is a cartesian closed category with a natural numbers object and finite sums, we can regard $\mathcal{C}$ as a tca as follows: the types will be the objects in $\mathcal{C}$, while $|X| = {\rm Hom}_\mathcal{C}(1, X)$ for any object $X$ in $\mathcal{C}$. In fact, this example would still work if we assumed that all the structure in $\mathcal{C}$ is \emph{weak} (by that we mean that we weaken the universal property by only requiring existence of a certain arrow; we drop the requirement that that arrow is also unique). But if $\mathcal{C}$ is genuinely cartesian closed and also \emph{well-pointed} (in that two parallel arrows $f, g: X \to Y$ will be equal whenever $fh = gh$ for any arrow $h: 1 \to X$), then the resulting tca will be extensional.
    \end{enumerate}
\end{exa}

\begin{lem}
    \label{Lemma: tca decidable}
    Every tca $\tca{T}$ contains an element $d \in |N \to N \to N|$ such that for all $a, b \in \mathbb{N}$ we have that 
    $$
        d \overline{a} \overline{b} = \begin{cases}
                    \overline{0}, & \text{ if } a \neq b, \\
                    \overline{1}  & \text{ if } a = b.
                \end{cases}
    $$
\end{lem}
\proof
    It is easy to check that
    $
        d := \mathsf{R}_{0 \to 0} (\mathsf{R}_0 (\mathsf{succ} \, \mathsf{0}) \mathsf{k}_{0,0}) (\lambda x y. \mathsf{R}_0 \mathsf{0} (\lambda z w.yz))
    $
    works. 
\qed

\begin{prop}
    Let $\tca{T}$ be a consistent tca. Then the following hold:
    \begin{enumerate}
        \item[(i)] The map $\mathbb{N} \to |N|: n \mapsto \overline{n}$ is injective.
        \item[(ii)] The maps $|A| \to |A + B|: a \mapsto \mathsf{inl} \, a$ and $|B| \to |A + B|: b \mapsto \mathsf{inr} \, b$ have disjoint images.
    \end{enumerate}  
\end{prop}
\proof
    (i) follows from the previous lemma. For (ii), consider \[ h := \mathsf{case}(\lambda x.\overline{0})(\lambda x.\overline{1}) \in |A + B \to N| \] and note that $h(\mathsf{inl} \, a) = \overline{0}$ while $h(\mathsf{inr} \, b) = \overline{1}$.
\qed

\subsection{Assemblies}

Following Longley \cite{longley99}, we generalise the usual category of assemblies over a pca, and define a category of assemblies on a tca $\tca{T}$.

\begin{defi}
    Let $\tca{T}$ be a tca. An \emph{assembly} on $\tca{T}$ is a triple $(X,A,\alpha)$, where $X$ is a set, $A$ a type of $\tca{T}$, and $\alpha: X \to \Pow^*(|A|)$ (with $\Pow^*(|A|)$ being the collection of inhabited subsets of $|A|$). A morphism of assemblies $f: (X,A,\alpha) \to (Y,B,\beta)$ is a function $f: X \to Y$ such that there exists an element $e \in |A \to B|$ with the property that whenever $a \in \alpha(x)$, then $ea \in \beta(f(x))$. We say that $e$ \emph{witnesses} or \emph{tracks} $f$, and refer to $(X,A,\alpha)$ with $X$.
\end{defi}

We denote the resulting category by $\Asm_\tca{T}$. We might omit any mention of the tca $\tca{T}$ whenever it is clear from the context or irrelevant which specific tca we are referring to. 

Our immediate goal is now to show that $\Asm_\tca{T}$ is a cartesian closed category with a natural numbers object and finite disjoint coproducts. We will content ourselves with describing all the relevant constructions, leaving a verification of their correctness to the reader. To this end, let $(X,A,\alpha)$ and $(Y,B,\beta)$ be assemblies.

\begin{description}
    \item[Initial Object] The initial object is $(\emptyset, \bot, \emptyset)$.
    \item[Terminal Object] The terminal object is $(\set{0}, \top, \phi)$, where $\phi$ is the map $0 \mapsto \set{\mathsf{t}}$.
    \item[Product] The product of $X$ and $Y$ is $(X \times Y, A \times B, \alpha \times \beta)$, where $\alpha \times \beta$ is the map $(x,y) \mapsto \set{k \, : \, \fst \, k \in \alpha(x), \snd \, k \in \beta(y)}$.
    \item[Pullback] Given morphisms $f: (X,A,\alpha) \to (Z,C,\gamma)$ and $g: (Y,B,\beta) \to (Z,C,\gamma)$, their pullback is $(P,A \times B, \chi \upharpoonright P)$, where 
    $$
        P = \Set{(x,y) \in X \times Y}{f(x) = g(y)},
    $$
    $\chi$ is as before and the maps $P \to X$ and $P \to Y$ are the projections. 
    \item[Coproduct] The coproduct of $X$ and $Y$ is $(X + Y, A + B, \gamma)$, where $\gamma$ is the map defined by $\inl(x) \mapsto \Set{\inl \, z}{z \in \alpha(x)}$ and $\inr(y) \mapsto \Set{\inr \, z}{z \in \beta(y)}$. Note that finite coproducts are stable under pullback.
    \item[Equaliser] Let $f, g: (X,A,\alpha) \to (Y,B,\beta)$. The equaliser of $f$ and $g$ is $(E, A, \alpha \upharpoonright E)$ where $E = \Set{x \in X}{f(x) = g(x)}$ and $e: E \to X$ is the inclusion.
    \item[Natural Numbers Object] Consider $(\mathbb{N}, N, \nu)$, where $\nu$ is the map $n \mapsto \set{\numeral{n}}$. The map $z: 1 \to \mathbb{N}$ is given by $0 \mapsto 0$ and witnessed by $\lambda x.\mathsf{0} \in |\top \to N|$.  The map $s: \mathbb{N} \to \mathbb{N}$ is defined by $n \mapsto n + 1$ and witnessed by $\mathsf{succ}$.
    \item[Exponential] We have that $X^Y = (Z, B \to A, \alpha^\beta)$, where $Z$ is the set of morphisms $Y \to X$ and $\alpha^\beta(f)$ is the set of elements tracking $f$.
\end{description}

\begin{thm}
    The category $\Asm_\tca{T}$ of assemblies is a cartesian closed category with a natural numbers object and finite disjoint coproducts. \qed
\end{thm}

\subsection{A hyperdoctrine} Our next goal is to show that $\Asm_\tca{T}$ is a Heyting category. To that purpose a simplified presentation of the subobject hyperdoctrine of $\Asm_\tca{T}$ will prove useful. 

\begin{defi}
    Let $(X, A, \alpha)$ be an assembly. A \emph{predicate on $(X, A, \alpha)$} is a tuple $(B,\beta)$ consisting of a type $B$ and a map $\beta: X \to \P(|B|)$ for which there is an element $f \in |B \to A|$ such that $b \in \beta(x)$ implies $f \cdot b \in \alpha(x)$.
    
    We say that a predicate $(B,\beta)$ is a \emph{subpredicate} of $(C,\gamma)$ and write $(B,\beta) \leq (C,\gamma)$ if there is an element $f \in |B \to C|$ such that for all $x \in X$, if $b \in \beta(x)$, then $f \cdot b \in \gamma(x)$. Note that this defines a \emph{preorder of predicates on $(X, A, \alpha)$}.
\end{defi}

\begin{prop}
    \label{Proposition: Equivalence of Categories}
    The preorder of subobjects of $(X, A, \alpha)$ as a thin category and the category of predicates of $(X, A, \alpha)$ are equivalent.
\end{prop}
\proof
    Given a predicate $({B}, \beta)$ on some assembly $(X, A, \alpha)$ as witnessed by $f_B \in |B \to A|$, we can define an assembly
    $
        (Y_B, {B}, \beta),
    $
    where $Y_B = \Set{x \in X}{\beta(x) \neq \emptyset}$.
    Note that $(Y_B, B, \beta)$ is a subobject of $(X, A, \alpha)$, witnessed by $f_B$. If $({B},\beta) \leq ({C},\gamma)$ witnessed by some $f$, then the inclusion $Y_B \subseteq {Y}_C$ is witnessed by $f$, i.e., ${Y}_B \subseteq {Y}_C$ as subobjects of $(X, A, \alpha)$.

    Conversely, if $(Y, B, \beta)$ is a subobject of $(X, A, \alpha)$, i.e., there is a monomorphism $f: (Y, B, \beta) \to (X, A, \alpha)$, then we can obtain a predicate $(C_Y, \gamma_Y)$ on $(X, A, \alpha)$, where $C_Y = B$ and $\gamma_Y(x) = \bigcup_{y \in f^{-1}(x)} \beta(y)$ for $x \in X$. If $(Y, B, \beta) \subseteq (Z, D, \delta)$ as subobjects of $(X,A, \alpha)$, then there is an element $f \in  |C \to D|$ witnessing this fact. By our definitions, this morphism also witnesses that $(C_Y, \gamma_Y) \leq (C_Z, \gamma_Z)$.

    In fact, we just defined functors ${Y}_{(-)}$ and $B_{(-)}$ between the thin categories $\Sub(X, A, \alpha)$ of subobjects of $(X, A, \alpha)$ and predicates of $(X, A, \alpha)$. A routine check shows that $B_{{Y}_{(B,\beta)}} = (B, \beta)$ and ${Y}_{B_{(Y,B,\beta)}} \cong {Y}$. We can thus conclude that both functors are full, faithful and essentially surjective, i.e., an equivalence of categories.
\qed

Let us start by noting that the predicates on $(X, A, \alpha)$ form a pre-Heyting algebra. From the previous proposition it then follows that $\Sub(X, A, \alpha)$ is a pre-Heyting algebra as well.
\begin{description}
    \item[Top Element] The top element of $P(X, A, \alpha)$ is $(A,\alpha)$.
    \item[Bottom Element] The bottom element of $P(X, A, \alpha)$ is $(\bot,\lambda x.\emptyset)$.
    \item[Conjunction] We have that $(B, \beta) \wedge (C, \gamma) = (B \times C, \beta \times \gamma)$.
    \item[Disjunction] $(B, \beta) \lor (C,\gamma) = (B + C, \delta)$ with $\delta(x) = \{ \inl \, y \, : \, y \in \beta(x) \} \cup \{ \inr \, y \, : \, y \in \gamma(x) \}$.
    \item[Implication] $(C, \gamma) \Rightarrow (B, \beta) = (A \times (C \to B), \delta)$, where $m \in \delta(x)$ whenever $\fst \, m \in \alpha(x)$ and $\snd \, m \in |C \to {B}|$ such that if $n \in \gamma(x)$, then $(\snd \, m) \, n \in \beta(x)$. 
\end{description}

We are now ready to define a hyperdoctrine $P: \Asm_\tca{T}^{\mathrm{op}} \to \mathbf{preHA}$ (where $\mathbf{preHA}$ is the category of pre-Heyting algebras) such that $P(X, A, \alpha)$ is the collection of all predicates on $(X, A, \alpha)$.

If $f: (Y, B, \beta) \to (X, A, \alpha)$ is a morphism of assemblies, then we define \[ Pf: P(X, A, \alpha) \to P(Y, B, \beta) \] by stipulating that:
$$
    (Pf)(C, \gamma) = (C \times B, \gamma_f),
$$
where for $y \in Y$ we define that:
$$
    \gamma_f(y) = \gamma(f(y)) \times \beta(y) := \{ \,  \mathsf{pair} \, m \, n \, : \, m \in \gamma(f(y)), n \in \beta(y) \, \}
$$
Note that $(Pf)(C,\gamma) = B_{f^* Y_{(C, \gamma)}}$ if $f^*Y_{(C, \gamma)}$ is the pullback of the subobject $Y_{(C, \gamma)}$ along $f$.

We define $\exists_f: P(Y, B, \beta) \to P(X, A, \alpha)$ by $$\exists_f(C,\gamma) = (C,\gamma_\exists),$$ where $$\gamma_\exists(x) = \bigcup_{y \in f^{-1}(x)} \gamma(y).$$

\begin{prop}
    The morphism $\exists_f$ is the left adjoint of $P(f)$.
\end{prop}
\proof
    As we are working with thin categories, it suffices to show that the following equivalence holds:
    $$
        \exists_f (C, \gamma) \leq_{P(X, A, \alpha)} (D, \delta) \iff (C, \gamma) \leq_{P(Y, B, \beta)} (Pf)(D, \delta)
    $$
    for $(C,\gamma) \in P(X, A, \alpha)$ and $(D,\delta) \in P(Y, B, \beta)$.
    
    For the first direction, assume that $\exists_f (C, \gamma) \leq (D, \delta)$. This means that there is an element $g \in |C \to D|$ such that $n \in \gamma_\exists(x)$ implies that $g \cdot n \in \delta(x)$. To show that $(C, \gamma) \leq (Pf)(D, \delta)$ it suffices to provide an element $h \in | C \to D \times B|$ such that $m \in \gamma(y)$ implies $hm \in \delta_f(y)$. Let $e_C \in | C \to B|$ be a witness for $(C, \gamma) \in P(Y, B, \beta)$, and take $h := \lambda m. \pair (gm)(e_C m)$. It follows that if $m \in \gamma(y)$, then $hm \in \delta(f(y)) \times \beta(y) = \beta_f(y)$. 
    
    For the other direction, assume that $(C,\gamma) \leq (Pf)(D,\delta)$, i.e., there is an element $h \in | C \to D \times B|$ such that if $n \in \gamma(y)$, then $hn \in \delta(f(y)) \times \beta(y)$. Now, if $n \in \gamma_\exists(x)$, then there is some $y \in f^{-1}(x)$ such that $n \in \gamma(y)$. Hence, $\fst \, (h \, n) \in \delta(f(y)) = \delta(x)$. This shows that $g := \lambda n. \fst \, (h \, n)$ witnesses that $\exists_f (C, \gamma) \leq (D, \delta)$.
\qed

The map $\forall_f: P(Y,B,\beta) \to P(X,A,\alpha)$ is defined by $\forall_f(D,\delta) = (A \times (B \to D), \delta_\forall)$, where $m \in \delta_\forall(x)$ if and only if $\fst \, m \in \alpha(x)$ and $\snd \, m \in |B \to D|$ has the property that $(\snd \, m) \,  n \in \delta(y)$ for all $n \in \beta(y)$ and $y \in f^{-1}(x)$.

\begin{prop}
    The morphism $\forall_f$ is the right adjoint of $P(f)$.
\end{prop}
\proof
    It suffices to show that the following equivalence holds:
    $$
        (C, \gamma) \leq_{P(X, A, \alpha)} \forall_f(D, \delta) \iff (Pf)(C, \gamma) \leq_{P(Y, B, \beta)} (D, \delta)
    $$
    for $(C, \gamma) \in P(X)$ and $(D, \delta) \in P(Y)$.
    
    For the first direction, suppose that $(C, \gamma) \leq_{P(X,A,\alpha)} \forall_f(D,\delta)$ is witnessed by $g \in | C \to A \times (B \to D)|$ such that if $n \in \gamma(x)$, then $gn \in \delta_\forall(x)$.  Let $h := \lambda m. (\snd \, g (\fst \, m)) (\snd \, m)$. If $m \in \gamma_f(y)$, then $\fst \, m \in \gamma(f(y))$ and $\snd \, m \in \beta(y)$. Hence, $g (\fst \, m) \in \delta_\forall(f(y))$ such that $h m = (\snd \, (g (\fst \, m)))(\snd \, m) \in \delta(y)$.
    
    For the second direction, assume that $(Pf)(C, \gamma) \leq (D, \delta)$. This means that there is an element $g \in | C \times B \to D|$ such that if $k \in \gamma_f(y)$, then $gk \in \delta(y)$. We show that $h := \lambda k. \pair \, (e_C k) \, (\lambda l. g (\pair \, k \, l))$ witnesses that $(C, \gamma) \leq \forall_f(D,\delta)$: Let $k \in \gamma(x)$. As $(C, \gamma) \in P(X, A, \alpha)$, let $e_C \in | C \to B|$ be a witness for $(C, \gamma) \in P(Y, B, \beta)$. It follows that $e_C k \in \alpha(x)$. Furthermore, let $l \in \beta(y)$ for some $y \in f^{-1}(x)$. Then $\pair \, k \, l \in \gamma(f(y)) \times \alpha(x) = \gamma_f(y)$, and hence, $g (\pair \, k \, l) \in \delta(y)$. These calculations show that $h k \in \delta_\forall(x)$.
\qed

\begin{prop}
    The morphisms $P(f)$, $\exists_f$ and $\forall_f$ satisfy the Beck-Chevalley condition, i.e., if
    \begin{center}
        \begin{tikzcd}
            X \arrow[r, "f"] \arrow[d, "g"] & Y \arrow[d, "h"] \\
            Z \arrow[r, "k"] & W 
        \end{tikzcd}
    \end{center}
    is a pullback square in $\Asm$, then the squares
    \begin{center}
        \begin{tikzcd}
            P(X) \arrow[r, "\exists_f"] & P(Y)  \\
            P(Z) \arrow[u, "P(g)"] \arrow[r, "\exists_k"] & P(W) \arrow[u, "P(h)"]
        \end{tikzcd}
        \hspace{5mm} and \hspace{5mm}
        \begin{tikzcd}
            P(X) \arrow[r, "\forall_f"] & P(Y)  \\
            P(Z) \arrow[u, "P(g)"] \arrow[r, "\forall_k"] & P(W) \arrow[u, "P(h)"]
        \end{tikzcd}
    \end{center}
    both commute.
\end{prop}
\proof
    We leave the verification that the left hand square commutes to the reader; from this the commutativity of the right hand square follows by adjointness.
\qed

\begin{thm}
    The category $\Asm_\tca{T}$ of assemblies is a Heyting category. \qed
\end{thm}

We conclude this section with a lemma which will prove useful later:
\begin{lem} \label{Lemma: modified negation}
    Suppose $(C, \gamma)$ is a predicate on an assembly $(X, A, \alpha)$. If $|\bot|$ is inhabited, then we can we put $\lnot(C,\gamma) = (A, \beta)$ with $\beta(x) = \bigcup \{ \alpha(x) \, : \, \gamma(x) = \emptyset \}$.
\end{lem}
\proof
    Suppose $\mathsf{f} \in |\bot|$. Note that according to the definition above \[ (C, \gamma) \to (\bot,\lambda x.\emptyset) = (A \times (C \to \bot), \delta) \] with $m \in \delta(x)$ whenever $\fst \, m \in \alpha(x)$, $\snd \, m \in |C \to \bot|$ and $\gamma(x)$ is empty. This means that $(C, \gamma) \to (\bot,\lambda x.\emptyset) \leq (A, \beta)$ is witnessed by $\fst \in |A \times (C \to \bot) \to A|$, while $(A, \beta) \leq (C, \gamma) \to (\bot,\lambda x.\emptyset)$ is witnessed by $\lambda x. \mathsf{pair} \, x \, (\lambda y.\mathsf{f}) \in |A \to (A \times (C \to \bot))|$.
\qed

\begin{rem} The assumption that $|\bot|$ is inhabited may sound paradoxical. However, there are many examples of tcas which satisfy this assumption: for instance, the collection of inhabited subsets of a pca $\mathcal{A}$. Indeed, tcas for which $|\bot|$ is inhabited play an important role in this paper: the reason is that we are interested in forms of modified realizability. The idea behind the category-theoretic treatment of modified realizability is that there is a distinction between the actual and potential realizers of a statement, where every statement, including $\bot$, has a potential realizer; of course, $\bot$ does not have an actual realizer (see \cite{vanOosten2008}). In the categories of assemblies that the we will discuss below the elements of $|\bot|$ are the \emph{potential} realizers of $\bot$.
\end{rem}

\section{Assemblies for modified realizability}

In this section we will take a closer look at the category of assemblies $\Asm_\tca{T}$ over an \emph{extensional} tca $\tca{T}$. For some interesting examples of extensional tcas, we refer to \cite{longley99,LietzStreicher2002,longleynormann15} as well as the next section. We will show that in that case we can regard $\tca{T}$ as a category which embeds into $\Asm_\tca{T}$. In addition, in the internal logic of $\Asm_\tca{T}$ the characteristic principles of modified realizability hold: the axiom of choice for finite types as well as the independence of premise principle if $|\bot|$ is inhabited. This generalises the contents of Chapter 2 in the MSc thesis by Mees de Vries \cite{devries17} supervised by the first author.

\emph{Throughout this section $\tca{T}$ will be an extensional tca.}

\subsection{Embedding types into assemblies} From any extensional $\tca{T}$ we can construct a category which we will also denote by $\tca{T}$. The objects of this category are the types of $\tca{T}$ and its morphisms $A \to B$ are the elements of $|A \to B|$; we will also refer to these morphisms as \emph{morphisms of types}. The identity arrow is given by $\lambda x.x$ and the composition of $f$ and $g$ by $\lambda x.g(f(x))$; the axioms for a category follow from the extensionality of the tca $\tca{T}$. In fact, using extensionality one can show that $\tca{T}$ is a well-pointed and cartesian closed category, with products given by $A \times B$ and exponentials by $A \to B$. For that reason, we may also write $B^A$ instead of $A \to B$.

In addition, there is a functor
\[ E: \tca{T} \to \Asm_\tca{T} \]
which sends a type $A$ to $E(A) = (|A|, A, \alpha)$ with $\alpha(x) = \{ x \}$, while $(Ef)(x) = f(x)$. This functor $E$ is clearly full and faithful. We now take a closer look at its image. 

\begin{defi}
    An assembly $(X,A,\alpha)$ is called \emph{modest} if, for all $x,y \in X$, if $a \in \alpha(x)$ and $b \in \alpha(y)$ are such that $a = b$, then $x = y$. We say that $(X,A,\alpha)$ is \emph{strongly modest} if, additionally, for all $x \in X$ and $a,b \in \alpha(x)$ we have that $a = b$. We say that an assembly $(X,A,\alpha)$ is \emph{exhaustive} if for all $a \in |A|$, there is $x \in X$ such that $a \in \alpha(x)$. Finally, we say that an assembly is a \emph{base} or \emph{basic} if it is both exhaustive and strongly modest.
\end{defi}

\begin{thm}
    The functor $E: \tca{T} \to \Asm_\tca{T}$ is full and faithful. Moreover, all the assemblies in the image of this functor are basic and every basic assembly is isomorphic to an assembly in the image of this functor. 
\end{thm}
\proof
    Note that by definition $E(A)$ is a base, and a straightforward calculation shows that $E A \cong (X, A, \alpha)$ for any basic assembly $(X,A,\alpha)$.
\qed

In fact, more is true: the functor $E$ preserves the cartesian closed structure. To show this, we need the following lemma.

\begin{lem}[Lifting of Morphisms]
    \label{Lemma: Lifting of Morphisms}
    Let $(X, A, \alpha)$ be a strongly modest assembly and $(Y, B, \beta)$ be a modest exhaustive assembly. If $f: A \to B$ is a morphism of types, then there is a unique morphism $\bar f: (X, A, \alpha) \to (Y, B, \beta)$ of assemblies tracked by $f$.
\end{lem}
\proof
    We construct $\bar f$ as follows. Given $x \in X$, let $\bar{f}(x) \in Y$ be the unique element such that $f(a) \in \beta(\bar{f}(x))$ for all $a \in \alpha(x)$. The map $\bar f$ is well-defined because $(X, A, \alpha)$ is strongly modest and $(Y, B, \beta)$ is both modest and exhaustive. By construction $\bar f$ is witnessed by $f$ and hence, $\bar f$ is a morphism of assemblies.
    
    For the uniqueness, let $\bar f': X \to Y$ be a morphism of assemblies witnessed by $f$. Let $x \in X$ and $a \in \alpha(x)$. Then $f(a) \in \beta(\bar f'(x))$ and $f(a) \in \beta(\bar f(x))$. By modesty of $Y$, it follows that $\bar f'(x) = \bar f(x)$. 
\qed

\begin{prop}
    In $\Asm_\tca{T}$ the assemblies which are modest are closed under finite products and exponentials, as are the basic assemblies. Moreover, $E$ preserves the cartesian closed structure.
\end{prop}
\proof
    First of all, let us note that $E(\top)$ is a terminal object in $\Asm_\tca{T}$ because $\mathsf{t}$ is the sole element of $|\top|$.
    
    Second, let $(X, A, \alpha)$ and $(Y, B, \beta)$ be two modest assemblies and consider their product $(X \times Y, A \times B, \alpha \times \beta)$. If both $(X, A, \alpha)$ and $(Y, B, \beta)$ are exhaustive, then $\bigcup \Im(\alpha) = |A|$ and $\bigcup \Im(\beta) = |B|$, where we have written $\Im(f)$ for the image of a map $f$. Then it follows that $\bigcup \Im(\alpha \times \beta) = |A| \times |B|$, i.e. $X \times Y$ is exhaustive. 
    
    For modesty, let $(x_0,y_0), (x_1,y_1) \in X \times Y$ and $(a_0, b_0) \in (\alpha \times \beta)(x_0,y_0)$ and $(a_1, b_1) \in (\alpha \times \beta)(x_1,y_1)$ be such that $(a_0,b_0) = (a_1,b_1)$. This means that $a_0 = a_1$ and $b_0 = b_1$ as well as $a_i \in \alpha(x_i)$ and $b_i \in \beta(y_i)$. So if $X$ and $Y$ are modest, it follows that $x_0 = x_1$ and $y_0 = y_1$. Hence $X \times Y$ is modest as well.

    For strong modesty, let $(x,y) \in X \times Y$ be such that $(a_0,b_0) \in (\alpha \times \beta)(x,y)$ and $(a_1,b_1) \in (\alpha \times \beta)(x,y)$. By our assumption that $X$ and $Y$ are strongly modest, it follows that $a_0 = a_1$, $b_0 = b_1$, and hence that $(a_0,b_0) = (a_1,b_1)$.
    
    Third, we consider exponentials. Again, let $(X,A,\alpha_A)$ and $(Y,B,\alpha_B)$ be two modest assemblies and let $(C, B \to A, \gamma)$ be their exponential, where $C$ is the set of morphisms $Y \to X$ and $\gamma(f)$ is the set of realizers for $f$. 
    
    For modesty of $X^Y$, let $c_0 \in \gamma(f)$ and $c_1 \in \gamma(g)$ for some $f, g \in C$ such that $c_0 = c_1$. Given $y \in Y$ and $b \in \beta(y)$, we have that $c_0(b) = c_1(b)$ and hence, by modesty of $X$, that $f(y) = g(y)$. It follows that $f = g$ and $X^Y$ is modest.
        
    Now suppose both $X$ and $Y$ are bases. To see that $X^Y$ is strongly modest, let $c_0,c_1 \in \gamma(f)$. Given any $b \in B$, find $y \in Y$ such that $b \in \beta(y)$. It follows that $c_0(b), c_1(b) \in \alpha(f(y))$. Hence $c_0(b) = c_1(b)$ and therefore $c_0 = c_1$ by extensionality. 
    
    Finally, by our assumptions and Lemma \ref{Lemma: Lifting of Morphisms}, we know that $\bigcup \Im(\gamma) = |B \to A|$, and so, $X^Y$ is also exhaustive.
\qed

\begin{prop}
    If $\tca{T}$ is standard, then $N$ is a natural numbers object in $\tca{T}$. Moreover, this natural numbers object is preserved by $E$.
\end{prop}
\proof
  The natural numbers object in $\Asm$ is $(\mathbb{N},N ,\nu)$ with $\nu(n) = \set{\numeral{n}}$, which coincides with $E(N)$ if $\tca{T}$ is standard. And because $E$ is a full and faithful functor preserving the terminal object, it reflects the natural numbers object. So $N$ is the natural numbers object in $\tca{T}$.
\qed

\subsection{IP, AC and MP}
In this section, we investigate which common principles hold in the assemblies $\Asm_\tca{T}$ for a given extensional tca $\tca{T}$. 

\begin{thm}
    \label{Theorem: IP} 
    If $\tca{T}$ is an extensional tca for which $|\bot|$ is inhabited, then in the category of assemblies $\Asm_\tca{T}$ over that tca, the following independence of premise principle is satisfied:
    $$
        (\neg \phi \rightarrow \exists y^Y \ \psi(y)) \rightarrow (\exists y^Y \ \neg \phi \rightarrow \psi(y)),
    $$
    for all modest and exhaustive $Y$.
\end{thm}
\proof
    Let $(X, A, \alpha)$ be the assembly that is the context of the independence of premise principle and $(Y, B, \beta)$ be an assembly which is both modest and exhaustive. We will first calculate the predicates on $(X, A, \alpha)$ that correspond to the premise and conclusion of IP. To fix notation, let $(C, \gamma) = \heyting{ \phi} \in P(X)$ be the predicate corresponding to ${\phi}$, and $(D, \delta) = \heyting{\psi} \in P(X \times Y)$ be the predicate corresponding to $\psi$. Then we may assume that the interpretation of $\lnot \phi$ will be $(A, \mu)$ with $\mu(x) = \bigcup \{ \, \alpha(x) \, : \, \gamma(x) = \emptyset \} = \{ n \, : \, n \in \alpha(x) \mbox{ and } \gamma(x) = \emptyset \}$, by Lemma \ref{Lemma: modified negation}.
    
    It is straightforward to compute that $\heyting{\neg \phi \rightarrow \exists y^Y \psi(y)} = (D^A \times A, \theta)$, where $(m,n) \in \theta(x)$ if and only if $n \in \alpha(x)$ and $m: A \to D$ is a morphism of types such that $i \in \mu(x)$ implies $m i \in \bigcup_{y \in Y} \delta(x,y)$. 
    
    On the other hand, $\heyting{\exists y \, ({\neg\phi} \rightarrow \psi(y))} = (D^{A \times B} \times (A \times B), \eta)$ such that $(k,l) \in \eta(x)$ if and only if there is some $y \in Y$ such that $l \in (\alpha \times \beta)(x,y)$ and $k j \in \delta(x,y)$ for any $j \in \mu(x) \times \beta(y)$.
    
    To show that the above principle of independence of premise holds, it suffices to show that $\heyting{\neg\phi \rightarrow \exists y^Y \psi(y)} \leq \heyting{\exists y (\neg\phi \rightarrow \psi(y))}$. To do so, we have to provide a morphism of types $f: D^{A} \times A \to D^{A \times B} \times (A \times B)$ such that $(m,n) \in \theta(x)$ implies $f(m,n) \in \eta(x)$. 
    As $\heyting{\psi}$ is a predicate on $X \times Y$, there is a morphism of types $g: D \to A \times B$ such that $i \in \delta(x,y)$ implies $gi \in \alpha_{X \times Y}(x,y) = \alpha(x) \times \beta(y)$. Then take $f(m,n) := (\lambda  i. m n, (n,\snd (g (m n))))$. We have to show $(m,n) \in \theta(x)$ implies $f(m,n) \in \eta(x)$.
    
    So assume that $(m,n) \in \theta(x)$, i.e., $n \in \alpha(x)$ and $m: A \to D$ is a morphism of types such that $i \in \mu(x)$ implies $m i \in \bigcup_{y \in Y} \delta(x,y)$. We have $mn \in D$, $g(mn) \in A \times B$ and $\snd(g(mn)) \in B$. Since $Y$ is exhaustive, there is an $y \in Y$ with $\snd(g(mn)) \in \beta(y)$. So $(n,\snd(g(mn))) \in (\alpha \times \beta)(x,y)$.

    It remains to show that $mn \in \delta(x,y)$ whenever $j \in \mu(x) \times \beta(y)$. But if $\fst \, j \in \mu(x)$, then $\gamma(x) = \emptyset$ and $\mu(x) = \alpha(x)$. Therefore  $n \in \alpha(x)$ and $mn \in \delta(x,y')$ for some $y' \in Y$. We would like to have $y = y'$, but this follows because $\snd \, (g(mn)) \in \beta(y')$ and $Y$ is modest.
\qed

\begin{thm}
    \label{Theorem: AC}
    Let $\tca{T}$ be an extensional tca. The category $\Asm_\tca{T}$ satisfies the axiom of choice for all basic assemblies. In particular, if $\tca{T}$ is also standard, we have the axiom of choice in all finite types.
\end{thm}
\proof
    Let $(X, A, \alpha)$ and $(Y, B, \beta)$ be basic assemblies, and $(Z, C, \gamma)$ be an arbitrary assembly. We will prove the axiom of choice for $X$ and $Y$ in context $Z$:
    $$
        (\forall x^X \exists y^Y \ \phi(x,y,z)) \rightarrow (\exists f^{X \to Y} \forall x^X \ \phi(x,fx,z)).
    $$
    We will compute the predicates corresponding to the premise and conclusion of the statement and then prove that the former is a subpredicate of the latter. To do so, assume that $(D, \delta) = \heyting{\phi(x,y,z)} \in P(X \times Y \times Z)$ is witnessed by $\iota_D$. 
    
    It is then straightforward to compute that 
    $$
        \heyting{\forall x^X \exists y^Y \ \phi(x,y,z)} = (C \times D^{A \times C}, \eta),
    $$ where $(n,m) \in \eta(z)$ if and only if $n \in \gamma(z)$ and $m: A \times C \to D$ is a morphism of types such that for every $k \in (\alpha \times \gamma)(x,z)$, there is some $y \in Y$ with $mk \in \delta(x,y,z)$.
    Moreover, we compute that $$
        \heyting{\exists f^{X \to Y} \forall x^X \ \phi(x,fx,z)} = ((B^A \times C) \times (D \times 
        (A \times B^A \times C))^{A \times B^A \times C}, \theta),
    $$ 
    where $(n, m) \in \theta(z)$ if and only if there is some $f \in Y^X$ such that $n \in (\beta^\alpha \times \gamma)(f,z)$ and $m: A \times B^A \times C \to D \times (A \times B^A \times C)$ is a morphism of types such that for all $x \in X$ and $k \in (\alpha \times \beta^\alpha \times \gamma)(x,f,z)$ we have that $mk \in \delta(x,fx,z) \times (\alpha \times \beta^\alpha \times \gamma)(x,f,z)$.
    
    Let $z \in Z$ and $(n,m) \in \eta(z)$ be given, and define $\check f_{n,m} := \lambda j. \pi_B (\iota_D (m (j,n)))$. It follows that $\check f_{n,m}$ is a morphism of types from $A$ to $B$ as it is defined by a $\lambda$-expression. By Lemma \ref{Lemma: Lifting of Morphisms}, there is a morphism $f_{n,m}: X \to Y$ of assemblies which is tracked by $\check f$. Hence $(\check f, n) \in (\beta^\alpha \times \gamma)(f_{n,m}, z)$. 
    
    Furthermore, let $x \in X$ be given and assume that $k \in (\alpha \times \beta^\alpha \times \gamma)(x,f_{n,m},z)$. Then $m(\pi_A k, n) \in \delta(x, y, z)$ for some $y \in Y$. Now since 
    $$
        \check f_{n,m} (\pi_X k) = \pi_Y (\iota_{D} (m(\pi_X k, n))),
    $$ 
    it follows that $\check f_{n,m} (\pi_X k) \in \beta(y)$. By the fact that $Y$ is a base, it follows that $y = f_{n,m} x$. Hence, $m(\pi_X k, n) \in \delta(x,f_{n,m} x,z)$.
    
    The above reasoning shows that the following map does the job:
    $$
        \lambda (n,m) . ((\check f_{n,m},n), \lambda k. (m(\pi_X k, n), k) )
    $$
    This map is clearly a morphism of types as it is defined as a $\lambda$-expression.
\qed

Let $\tca{T}$ be a consistent tca. We say that an element $f$ of type $(N \times N) \to N$ \emph{solves the halting problem} if and only if the following condition holds for all natural numbers $a$ and $b$:
$$
    f(\numeral{a},\numeral{b}) = 
    \begin{cases}
        \Nsucc \, \mathsf{0}, & \text{ if the machine $a$ halts on input $b$,} \\
        \mathsf{0}, & \text{ otherwise.}
    \end{cases}
$$
For ease of notation we will write ``$ab {\downarrow}$'' to say that machine $a$ halts on input $b$, and ``$ab {\uparrow}$'' for the negation of the previous statement.

\begin{thm}
    \label{Theorem: MP fails}
    Let $\tca{T}$ be an extensional and standard tca that does not contain an element solving the halting problem. Then Markov's principle fails in $\Asm_{\tca{T}}$.
\end{thm}
\proof
    Let $(\mathbb{N}, N, \nu)$ be the natural numbers object. Because $\tca{T}$ is assumed to be standard, $\nu$ is an isomorphism and we need not distinguish between numerals $\numeral{a}$ and natural numbers $a$.
    
    Let $T = (N \times N \times N, \tau)$ be the predicate on the assembly $\mathbb{N} \times \mathbb{N} \times \mathbb{N}$ defined by
    $$
        \tau(a,b,n) = \{ \langle a, b, n \rangle \, : \, T(a,b,n) \} 
    $$
    where $T$ is Kleene's $T$-predicate. Suppose, for contradiction, that Markov's principle holds on $\Asm_\tca{T}$. It follows that $(A, \alpha) = \heyting{\neg\neg\exists n \, T(a,b,n)}$ is a subobject of $(B, \beta) = \heyting{\exists n \, T(a,b,n)}$. Note that $\beta(a,b) = \{ \langle \numeral{a},\numeral{b},\numeral{n} \rangle \, : \, T(a,b,n) \}$ and hence $\alpha(a,b) = \{ \langle \numeral{a}, \numeral{b}, \numeral{m} \rangle  \, : \, \lnot \, \lnot \exists n \in \mathbb{N} \, T(a,b,n) \}$ by Lemma \ref{Lemma: modified negation}. Therefore we are assuming that there is a morphism of types $g: N \times N \times N \to N \times N \times N$ such that
    \[ \lnot \lnot \exists n \, T(a,b,n) \to T(a,b,g(a,b,m)). \] 
    Therefore $\exists n \, T(a,b, n) \leftrightarrow T(a,b,g(a,b,0))$ and the halting problem can be decided in $\tca{T}$ by checking $T(a,b,g(a,b,0))$; in other words, the halting problem is decidable in $\tca{T}$.
\qed

\begin{rem} We will recover the results in Chapter 2 of the MSc thesis of Mees de Vries \cite{devries17} if we restrict attention to the extensional tca of inhabited PERs over $K_1$ (see also the next section). In the other chapters De Vries shows that the resulting category of assemblies arises as the $\lnot\lnot$-separated objects of a realizability topos which combines aspects of both the extensional and modified realizability topos (for which, see \cite{vanOosten2008}).
\end{rem}

\section{Apartness assemblies and converse extensionality}

In this section we will introduce a new extensional tca, which we will call the \emph{apartness types}; in fact, we will be able to define a tca of apartness types for each consistent tca. Our main goal is to show that converse extensionality principles will hold in the category of assemblies over the apartness types.

\emph{Throughout this section $\tca{T}$ will be a consistent typed combinatory algebra.}

\subsection{Apartness types} In this section, we define a category of apartness types based on a typed combinatory algebra. To motivate this definition, consider the following constructions of \emph{PERs} over a tca $\tca{T}$.

\begin{defi}
    Let $\tca{T}$ be a consistent tca. The \emph{category of PERs} over $\tca{T}$ is defined as follows: its objects are pairs $(A, R)$ where $A$ is a type in $\tca{T}$ and $R$ is a partial equivalence relation on $|A|$ (that is, a symmetric and transitive relation on $|A|$). The morphisms $(A, R) \to (B, S)$ are equivalence classes of elements $f \in |A \to B|$ such that $a_0 R a_1$ implies $fa_0 S fa_1$ for all $a_0, a_1 \in |A|$. Two such elements $f, g \in |A \to B|$ will be considered equivalent if $fa S ga$ for each $a \in |A|$ for which $aRa$ holds. The \emph{inhabited PERs} over $\tca{T}$ are those PERs $(A,R)$ for which there is at least one element $a \in |A|$ such that $aRa$.
\end{defi}    

One can check that the (inhabited) PERs form a well-pointed cartesian closed category with a natural numbers object and weak finite sums. The apartness types are a variation on this construction. The main idea is to replace the partial equivalence relations by a notion of apartness, which is to be thought of as a positive way of saying that two elements are inequivalent. And the fact that the notion of apartness is positive means in particular that to assert that two elements are apart requires evidence (it is ``proof relevant''). With this in mind, we have settled on the following definition.

\begin{defi}
    Let $\tca{T}$ be a consistent tca. A \emph{$\tca{T}$-apartness type} $\mathcal{A}$ is a tuple $(A,T_A,T^-_A,\apart_A)$ consisting of types $T_A, T_A^- \in \tca{T}$, a subset $A \subseteq |T_A|$ and a relation ${\apart_A} \subseteq A \times A \times |T_A^-|$. Instead of $(a_0,a_1,n) \in {\apart_A}$ we will write $n : a_0 \apart_A a_1$. We will require both $A$ and $|T_A^-|$ to be inhabited. In addition, some conditions have to be satisfied:
    \begin{enumerate}
        \item Reflexivity. There are no $a \in A$, $n \in A^-$ such that $n: a \apart_A a$.
        \item Symmetry. There is $s \in |T_A \to T_A \to T_A^- \to T_A^-|$ such that whenever $n: a_0 \apart_A a_1$, then $s a_0 a_1 n : a_1 \apart_A a_0$.
        \item Transitivity. There is $t \in |T_A \to T_A \to T_A \to T_A^- \to (T_A^- + T_A^-)|$ such that whenever $n : a_0 \apart_A a_1$ and $a_2 \in A$, then either $t a_0 a_1 a_2 n = \inl(m)$ and $m : a_0 \apart_A a_2$, or $t a_0 a_1 a_2 n = \inr(m)$ and $m : a_1 \apart_A a_2$.
    \end{enumerate}
\end{defi}

For every apartness type $\mathcal{A}$ we can define an equivalence relation $\sim_A$ on $A$ by saying:
$$
    a_0 \sim_A a_1 \text{ if and only if there is no $n$ such that } n : a_0 \apart_A a_1.
$$

\begin{defi}
    If $\mathcal{A} = (A, T_A, T_A^-, \apart_A)$ and $\mathcal{B} = (B, T_B, T_B^-, \apart_B)$ are $\tca{T}$-apartness types, then a \emph{premorphism from $\mathcal{A}$ to $\mathcal{B}$} is an element $f \in |(T_A \to T_B) \times (T_A \to T_A \to T_B^- \to T_A^-)|$ such that:
    \begin{enumerate}
        \item $(\fst \, f)(a) \in B$ for all $a \in A$, and
        \item for all $a_0, a_1 \in A$ and $n \in T_B^-$, if $n : (\fst \, f)(a_0) \apart_B (\fst \, f)(a_1)$, then $(\snd \, f) a_0 a_1 n : a_0 \apart_A a_1$. 
    \end{enumerate}
\end{defi}

If $f$ is a premorphism, we will usually write $f$ for $\fst \, f $ and $f^-$ for $\snd \, f$ without causing confusion. Note that for such a premorphism, the map $f: A \to B$ preserves the equivalence relation $\sim$ by condition (ii). We define an equivalence relation on premorphisms:
\begin{quote}
    $f \sim g$ if and only if $(\fst \, f)(a) \sim_B (\fst \, g)(a)$ for all $a \in A$.
\end{quote}

\begin{defi}
    The category $\ApType_\tca{T}$ of $\tca{T}$-apartness types is the category whose objects are $\tca{T}$-apartness types and whose morphisms are equivalence classes of $\tca{T}$-premorphisms.
\end{defi}

\begin{prop}
    The category $\ApType_\tca{T}$ of $\tca{T}$-apartness types is a well-pointed cartesian closed category with a natural numbers object, binary coproducts and a weak initial object.
\end{prop}
\proof
    Let $(A, T_A, T_A^-, \apart_A)$ and $(B, T_B, T_B^-, \apart_B)$ be $\tca{T}$-apartness types. Their product can be computed as follows:
    $$
        (A, T_A, T_A^-, \apart_A) \times (B, T_B, T_B^-, \apart_B) = (A \times B, T_A \times T_B, T_A^- + T_B^-, \apart_{A \times B}),
    $$
    where 
    \[ A \times B = \{ \, x \in |T_A \times T_B| \, : \, \mathsf{fst} \, x \in A, \mathsf{snd} \, x \in B \} \] 
    and $n : x_0 \apart_{A \times B} x_1$ if and only if either $n = \inl(m)$ with $m : \fst \, x_0 \apart_A \fst , x_1$ or $n = \inr(m)$ with $m : \snd \, x_0 \apart_B \snd \, x_1$.
    
    The exponential $(B, T_B, T_B^-, \apart_B) \to (A, T_A, T_A^-, \apart_A)$ is $(E,T_E,T_E^-,\apart_E)$, where
    \begin{align*}
        T_E &= (T_B \to T_A) \times (T_B \to T_B \to T_A^- \to T_B^-), \\
        E &= \{ \, \text{premorphisms } B \to A \, \}, \\
        T_E^- &= T_B \times T_A^-,
    \end{align*}
    and $n : f \apart_E g$ if and only if $\fst \, n \in B$ and $\snd \, n : f (\fst \, n) \apart_A g (\fst \, n)$. 
    
    The natural numbers object $\mathcal{N} = (\mathcal{N}, N, \top, \apart_N)$ is defined as follows:
    \begin{align*}
        \mathcal{N} &= \{ \, \overline{n} \, : \, n \in \mathbb{N} \} , \text{ and, } \\
        x : n &\apart_N m \text{ if and only if } n \neq m.
    \end{align*}
    Establishing that this is indeed the natural numbers object makes use of the primitive recursion operator $\mathsf{R}$ of the tca $\tca{T}$ at hand.
    
    The following is the coproduct $(A, T_A, T_A^-, \apart_A) + (B, T_B, T_B^-, \apart_B)$:
    $$
        (A + B, T_A + T_B, T_A^- \times T_B^-, \apart_{A + B}),
    $$
    where $A + B = \{ \mathsf{inl} \, a  \, : \, a \in A \} \cup \{ \mathsf{inr} \, b \, : \, b \in B \}$ and $n : \inl \, a_0 \apart_{A + B} \inl \, a_1$ if $\fst \, n : a_0 \apart_A a_1$, $n : \inr \, b_0 \apart_{A + B} \inr \, b_1$ if $\snd \, n : b_0 \apart_A b_1$, while $n : \inl \, a \apart_{A + B} \inr \, b$ and $n : \inr \, b \apart_{A + B} \inl \, a$ are always true.
    
    Finally, the apartness types have a terminal object $\mathbf{1}$, explicitly calculated as \[ (\set{\mathsf{t}},\top,\top,\emptyset); \] this object is also weakly initial.
\qed

\begin{cor}
    Let $\tca{T}$ be a consistent tca. The apartness types $\ApType_\tca{T}$ form an extensional and standard tca in which $|\bot|$ is inhabited.
\end{cor}
\proof
    This is immediate from the previous proposition and Example \ref{Example: tcas}(v). Note that $\ApType_\tca{T}$ is standard by construction and that $|\mathcal{A}| = A/{\sim}$ for any apartness type $\mathcal{A}$.
\qed

\begin{lem} \label{Lemma: haltingprobleminapptypes}
    If $\ApType_\tca{T}$ has a functional solving the halting problem, then so does $\tca{T}$.
\end{lem}
\proof
    If $f$ solves the halting problem in $\ApType_\tca{T}$, then $\fst \, f$ solves the halting problem in $\tca{T}$.
\qed

\begin{defi}
    Let $\tca{T}$ be a tca. The \emph{category $\ApAsm_\tca{T}$ of $\tca{T}$-apartness assemblies} is the category $\Asm_{\ApType_\tca{T}}$ of $\ApType_\tca{T}$-assemblies.
\end{defi}

\begin{rem}
    Note that the existence of a universal type on $\ApType_\tca{T}$ induces a universal type on $\tca{T}$. With \cite[Proposition 2.5]{LietzStreicher2002} it follows that $\ApType_{\mathcal{K}_1}$ and $\ApType_{\mathcal{K}_2^\mathrm{rec}}$ do not have universal types. By \cite[Theorem 4.2]{LietzStreicher2002}, this means that the ex/reg-completion of the apartness assemblies $\ApAss_{\mathcal{K}_1}$ and $\ApAss_{\mathcal{K}_2^\mathrm{rec}}$ are not toposes.
\end{rem}

\subsection{Converse extensionality}

Our next aim is to show that in the apartness assemblies certain \emph{converse extensionality principles} $\CE_n$ hold:
\begin{equation*}
    \exists X ^ {n + 2 \to n + 1 \to n + 1 \to n}  \, \forall \Phi^{n+2} \, \forall f^{n+1}, g^{n+1} \, (\Phi f \neq_0 \Phi g \rightarrow f (X \Phi f g) \neq_0 g (X \Phi f g)).
\end{equation*}
Note that because $\ApType_\tca{T}$ is an extensional and standard tca, the axiom of choice for all finite types holds in $\ApAsm_\tca{T}$ (see Theorem \ref{Theorem: AC}). Therefore the principle $\CE_n$ is equivalent to:
\begin{equation*}
    \forall \Phi^{n+2} \, \forall f^{n+1}, g^{n+1} \, \exists x^n \ (\Phi f \neq_0 \Phi g \rightarrow fx \neq_0 gx).
\end{equation*}

\begin{thm}
    \label{Theorem: CE_0}
    If $\tca{T}$ is a consistent tca, then
    converse extensionality $\CE_0$ holds in $\ApAsm_\tca{T}$.
\end{thm}
\proof
    Given a finite type $\tau$, we will write $(X_\tau, A_\tau, \alpha_\tau)$ for $\heyting{\tau}$. Note that all finite types are basic, so we may assume that $X_\tau = |A_\tau|$ and $\alpha_\tau(x) = \{ x \}$.

    We compute 
    $$
        \heyting{\Phi f \neq_0 \Phi g} = (A_2 \times A_1 \times A_1 \times A_0, \alpha),
    $$
    where
    $$
        \alpha(\Phi,f,g,x) = 
            \begin{cases}
                \{ (\Phi,f,g,x) \},
                &\text{ if } \Phi f \neq \Phi g,  \\
                \emptyset, 
                & \text{ otherwise,}
            \end{cases}
    $$ 
    and 
    $$
        \heyting{fx \neq_0 gx } = (A_2 \times A_1 \times A_1 \times A_0, \beta),
    $$ 
    where 
    $$
        \beta(\Phi,f,g,x) =
        \begin{cases}
            \{ (\Phi,f,g,x) \} ,
            & \text{ if } fx \neq gx, \\
            \emptyset, 
            & \text{ otherwise.}
        \end{cases}
    $$
    If we follow the definitions, then
    \begin{align*}
        & \heyting{\exists x^0 \ (\Phi f \neq_0 \Phi g \rightarrow f x \neq_0 g x } \\
        & \hspace{1cm} = 
        (
        (A_2 \times A_1 \times A_1 \times A_0)^{(A_2 \times A_1 \times A_1 \times A_0)} \times (A_2 \times A_1 \times A_1 \times A_0),
        \delta
        ),
    \end{align*}
    where $(m,k) \in \delta(\Phi,f,g)$ if and only if there is some $x \in \mathbb{N}$ such that 
    $k \in (\alpha_2 \times \alpha_1 \times \alpha_1 \times \alpha_0)(\Phi,f,g,x)$
    and 
    $
    m: A_2 \times A_1 \times A_1 \times A_0 \to A_2 \times A_1 \times A_1 \times A_0
    $ 
    is a morphism of types such that if $j \in \alpha(\Phi,f,g,x)$ then $mj \in \beta(\Phi,f,g,x)$. However, since $m$ has to act like the identity, we can simplify this to:
    \begin{align*}
        & \heyting{\exists x^0 \ (\Phi f \neq_0 \Phi g \rightarrow f x \neq_0 g x } \\
        & \hspace{1cm} = 
        (A_2 \times A_1 \times A_1 \times A_0,
        \gamma
        ),
    \end{align*}
    where $k \in \gamma(\Phi,f,g)$ if there is some $x \in \mathbb{N}$ such that $k = (\Phi,f,g,x)$ and $\Phi f \not= \Phi g$ implies $fx \not= gx$.

    We want to show that 
        $$
            \heyting{\exists x^0 \ (\Phi f \neq_0 \Phi g \rightarrow f x \neq_0 gx)}
        $$
    is the maximal element in $P(\heyting{2} \times \heyting{1} \times \heyting{1})$. To this end, we will construct a morphism $F: A_{2 \times 1 \times 1} \to A_0$ such that $(\Phi,f,g,F(\Phi, f,g) \in \gamma(\Phi,f,g)$.
    
    We let $F(\Phi, f, g) = x$, where $ x$ is the numeral computed using the following informally described algorithm: Let $( \Phi, f, g) \in X_{2 \times 1 \times 1}$. By Lemma \ref{Lemma: tca decidable}, we can check whether $\Phi f = \Phi g$. If this is the case, return $0$. If this is not the case, compute the least natural number $x \leq \fst \, (\Phi^-( f, g, \mathsf{t}))$ such that $f x \neq g x$, and return $x$. To understand this, observe that $\Phi^-(f,g, \mathsf{t})$ is a witness for the fact that $f$ and $g$ are apart. Therefore it follows from the construction of exponentials in the category of apartness types that $\fst \, (\Phi^-(f, g, \mathsf{t}))$ is an element $y$ such that $ f y \not= g y$, and $x$ is the smallest such. (Note that $x$ depends on $ f$ and $ g$ only, not on $\Phi^-$.) Therefore $(\Phi,f,g,F(\Phi,f,g)) \in \gamma(\Phi,f,g)$, as desired. 
    
    It remains to show that $F$ is a morphism of types; that is, we have to show that there is a map $F^-$ that witnesses reflection of apartness, i.e., if $m: F(\Phi_0,f_0,g_0) \apart F(\Phi_1, f_1, g_1)$ then $$F^-(({\Phi}_0, f_0, g_0),({\Phi}_1, f_1,g_1),m) : ({\Phi}_0,f_0,g_0) \apart ({\Phi}_1, f_1,g_1).$$
    
    We will again describe the algorithm $F^-$ informally. Writing $x_0 = F(\Phi_0,f_0,g_0)$ and $x_1 = F(\Phi_1,f_1,g_1)$, we know ${x}_0 \not= {x}_1$. It now suffices to compute values $j \in \{0, 1\}$ and $y$ such that $f_0 y \not= f_1 y$ if $j = 0$ and $g_0 y \not= g_1 y$ if $j = 1$. Let us say that $x_i$ is the smallest value of $\{ x_0,x_1 \}$ different from 0. Then $f_i x_i \not= g_i x_i$. If $x_{1-i}$ is also different from 0 then $x_i < x_{1-i}$ and $f_{1-i} x_i = g_{1-i} x_i$. Therefore either $f_0 x_i \not= f_1 x_i$ or $g_0 x_i \not= g_1 x_i$. If $x_{1-i} = 0$, then $f_{1-i} 0 \not= g_{1-i} 0$ or $f_{1-i}n = g_{1-i}n$ for all $n$. In the former case we have by our choice of $i$ that $f_i 0 = g_i 0$ and therefore either $f_0 0 \not= f_1 0$ or $g_0 0 \not = g_1 0$; in the latter case we have $f_{1-i}x_i = g_{1-i}x_i$ and therefore either $f_0 x_i \not= f_1 x_i$ or $g_0 x_i \not= g_1 x_i$.
  \qed

\begin{rem} In particular, in the apartness assemblies over G\"odel's $T$ (considered as a tca) the converse extensionality principle $\CE_0$ holds. This should be contrasted with Howard's negative result \cite{Troelstra344} saying that there is no term in G\"odel's $T$ witnessing $\CE_0$ directly. The reason for this discrepancy is that in the apartness types any type 2 functional $\Phi$ comes with a functional $\Phi^-$ witnessing the fact that $\Phi$ reflects apartness. This means that as soon as we discover that $\Phi f \not= \Phi g$ we can find a natural number $y$ such that $fy \not= gy$; in particular, we can find the least such $y$. Clearly, this $y$ could also have been found using unbounded search, would this be available. But having one such $y$ available means that the least such $y$ can be found using bounded search; as a result, the least such $y$ can be found in the apartness types over any tca, including those tcas where unbounded search is not available, like the term model of G\"odel's $T$.
\end{rem}

The same technique cannot be straightforwardly used to prove the converse extensionality principle for higher types. The reason is that, despite the fact that every function comes with a witness for the reflection of apartness that we have access to, the functional $X$ will have to be extensional, so its output cannot depend on which specific witness it was given as part of its input. With certain tca's, it is however possible to succeed with a slight modification.

\begin{defi}
    We say that $\tca{T}$ has a \emph{modulus of continuity} if there is a functional $M$ of type $2 \to 1 \to 0$ in $\tca{T}$ such that for every $f^2$ and $x^1, y^1$, we have that if $x \upharpoonright (Mfx) = y \upharpoonright (Mfx)$, then $fx = fy$. (Here $x \upharpoonright n$ refers to the finite sequence consisting of the first $n$ elements of $x$; that is, $\langle x(0), x(1), \ldots, x(n-1) \rangle$.)
\end{defi}

\begin{thm}
    \label{Theorem: CE_1}
    If $\tca{T}$ is a consistent tca with a modulus of continuity, then converse extensionality $\CE_1$ holds in $\ApAsm_\tca{T}$.
\end{thm}
\proof   
    The proof is essentially the same as the proof of the previous theorem with a slight modification in the algorithm needed for accommodating the higher types we deal with here. For that reason, we will only state the modified algorithm and leave the details to the reader.

    Let $M$ be the modulus of continuity of $\mathcal{T}$. We will also fix a computable bijection $h: \mathbb{N} \to \mathbb{N} ^ {< \mathbb{N}}$ between natural numbers and finite sequences of natural numbers. If $s$ is such a finite sequence of natural numbers, we will put $s^*: N \to N$ where
    $$
        s^* n = \begin{cases}
                    s(n), & n < \mathsf{length}(s), \\
                    0, & \text{otherwise}. \\
                \end{cases}
    $$

    In addition, we apply the same notational convention as in the previous proof, i.e. $\heyting{\tau} = (X_\tau,A_\tau,\alpha_\tau)$ for any finite type $\tau$.
    Analogous to the situation in the previous proof, we need to provide a morphism of types
    $
    F: A_3 \times A_2 \times A_2 \to A_1
    $
    such that $x = F(\Phi,f,g)$ is such that $\Phi f \not=_0 \Phi g \to f x \not =_0 gx$.
    
    We now informally describe the algorithm that implements $F$ as follows. Let $(\Phi, f, g) \in X_3 \times X_2 \times X_2$ be given. By Lemma \ref{Lemma: tca decidable}, we can check whether $\Phi f = \Phi g$. If this is the case, return $\lambda n. 0$. If this is not the case, proceed as follows: we compute $b = \fst (\Phi^- (f, g, \mathsf{t}))$ and observe that $fb \not= gb$. If we write $m = \max(M f b, M g b)$ we have $fb = f(b \upharpoonright m)^*$ and $gb = g(b \upharpoonright m)^*$ and therefore $f(b \upharpoonright m)^* \not= g(b \upharpoonright m)^*$. We now search through all sequences $s$ of length $\leq m$ using the ordering induced by $h$ until $f s^* \neq g s^*$, in which case we output $s^*$.
\qed

We apply the results of this section to some specific tca's. 

\begin{cor}
    \begin{enumerate}
        \item Let $\mathcal{K}_1$ be the tca arising from Kleene's first model. Then, \[ \ApAsm_{\mathcal{K}_1} \vDash \IP + \AC + \CE_0 \] but $\ApAsm_{\mathcal{K}_1} \not \vDash \MP$. 
        \item Let $\mathcal{K}_2^\mathrm{rec}$ be the recursive submodel of Kleene's second algebra $\mathcal{K}_2$. Then, \[ \ApAsm_{\mathcal{K}_2^\mathrm{rec}} \vDash  \IP + \AC + \CE_0 + \CE_1 \] but $\ApAsm_{\mathcal{K}_2^\mathrm{rec}} \not \vDash \MP$.
    \end{enumerate}
\end{cor}
\proof
    Note that both $\mathcal{K}_1$ and $\mathcal{K}_2^\mathrm{rec}$ do not contain an element solving the halting problem: both tca's contain only recursive functions. Furthermore, note that  $\mathcal{K}_2^\mathrm{rec}$ has a modulus of continuity by \cite[Theorem 2.6.3]{Troelstra344}; note that the functional constructed there for $\mathcal{K}_2$ also exists in $\mathcal{K}_2^\mathrm{rec}$ as it is a computable function. Then apply Theorems \ref{Theorem: IP}, \ref{Theorem: AC}, \ref{Theorem: MP fails}, \ref{Theorem: CE_0} and \ref{Theorem: CE_1} and Lemma \ref{Lemma: haltingprobleminapptypes}.
\qed

\begin{rem} Note that $\eha + \AC + \MP$ proves $\CE_0$. Indeed, if $\MP_n$ is Markov's principle for objects of type $n$, then $\eha + \AC + \MP_n \vdash \CE_n$, as the following calculation shows:
    \begin{eqnarray*}
        \forall \Phi^{n+2} \, \forall f, g \, \big( \,  f =_0 g \to \Phi f =_0 \Phi g \, \big) & \xRightarrow{\EXT} \\
        \forall \Phi^{n+2} \, \forall f, g \, \big( \, \forall x^n \, f x  =_0 g x \to \Phi f =_0 \Phi g \, \big) & \xRightarrow{\IL}\\
        \forall \Phi^{n+2} \, \forall f, g \, \big( \Phi f \not=_0 \Phi g \to \lnot \forall x^n f x  =_0 g x \, \big) & \xRightarrow{\MP_n} \\
        \forall \Phi^{n+2} \, \forall f, g \, \big( \, \Phi f \not=_0 \Phi g \to \exists x^n \, f x  \not=_0 g x  \, \big) & \xRightarrow{\IL} \\
        \forall \Phi^{n+2} \, \forall f, g \, \exists x \, ( \Phi f \not=_0 \Phi g \to f x  \not=_0 g x  \, ) & \xRightarrow{\AC} \\
        \exists X \, \forall \Phi^{n+2} \, \forall f, g \, \big( \, \Phi f \not=_0 \Phi g \to f (X\Phi f g)  \not=_0 g (X \Phi f g)  \, \big), 
        \end{eqnarray*}
    where $\EXT$ stands for the axiom of extensionality and $\IL$ for intuitionistic logic.
    So what the previous corollary shows is that the implication cannot be reversed, in that $\eha + \AC + \CE_n \not\vdash \MP_n$ for $n = 0,1$. 
\end{rem}    

\begin{rem} Reasoning classically, a stronger result is possible, as pointed out to us by Ulrich Kohlenbach. Using modified realizability in ECF, we obtain a model of $\eha + \IP + \AC + \CE_n + \lnot\MP$. Under this interpretation, Markov's principle fails because witnessing 
    \[ \forall f^1 \, (\neg\neg \exists n^0 \, f(n)=0  \to \exists n^0 \, f(n)=0) \]
requires a discontinous functional. In addition, in ECF there are functionals witnessing $\CE_n$ (see \cite[Appendix]{Troelstra344}).
\end{rem}

\section{Conclusion and directions for future research}

We set out to witness the converse extensionality principles $\CE_n$ using Brouwer's notion of apartness. We have succeeded in witnessing $\CE_0$ using a modified realizability interpretation using apartness types. The interesting aspect of this fact is that we can find this witness in G\"odel's $T$, despite Howard's result saying that $\CE_0$ cannot be directly witnessed in the term model. Witnessing higher $\CE_n$ proved difficult, because we have to make sure that the witnessing functional $X$ is still an extensional function. One natural direction, then, might be to drop the requirement that the witnessing functional $X$ be extensional. 

As suggested to us by Ulrich Kohlenbach, it may also be interesting to use enrichment of data to witness restricted forms of extensionality. For instance, in proof mining one often needs extensionality in the form 
\[ x=_X p \land T(p)=_X p \to T(x)=_X x, \]
where $(X,d)$ is some metric space. Even when there are no functionals witnessing full extensionality, functionals witnessing such restricted forms of extensionality might exist (see, for instance, \cite{kohlenbachetal18}).

Finally, another open question is whether the apartness assemblies arise as a subcategory of a suitable realizability topos.

\bibliographystyle{alphaurl}
\bibliography{bibliography}

\appendix

\section{A proof-theoretic account}

Some of the results in this paper can also be stated in proof-theoretic language. We include those statements here for the benefit of those readers who are primarily interested in the proof-theoretic aspects of our work.

Let us work in $\eha$ where we include both binary sum and product types in the type structure. It will be convenient to work with a version in which equality at all finite types is primitive and we have combinators as on page 3 of our paper (hence $\lambda$-abstraction is implemented using the combinators ${\bf k}$ and ${\bf s}$).

By induction on this type structure we define for each finite type $\sigma$ two finite types $\sigma^+$ and $\sigma^-$, as well as two formulas ${\rm dom}_\sigma(x)$ and ${\rm app}_\sigma(x,y,z)$ in the language of $\eha$ where ${\rm dom}_\sigma(x)$ has one free variable $x$ of type $\sigma^+$ and ${\rm app}_\sigma(x,y,z)$ has free variables $x,y,z$ of types $\sigma^+, \sigma^+$ and $\sigma^-$, respectively. (Here ${\rm dom}_\sigma(x)$ stands for the ``domain'' of $\sigma$: it applies to those elements of type $\sigma^+$ which are suitable for interpreting elements of type $\sigma$. The predicate ${\rm app}_\sigma(x,y,z)$ should be read as: $z$ is evidence for the statement that $x$ and $y$ are apart; the variable $z$ is of type $\sigma^-$, so $\sigma^-$ is the type of evidence of elements of type $\sigma$ being apart.)

\allowdisplaybreaks
\begin{eqnarray*}
N^+ & :\equiv & N \\
N^- & :\equiv & \top \\
{\rm dom}_N(x) & :\equiv & x =_0 x \\
{\rm app}_N(x,y,z) & :\equiv & x \not = y \\ \\
(\sigma \times \tau)^+ & :\equiv & \sigma^+ \times \tau^+ \\
(\sigma \times \tau)^- &:\equiv & \sigma^- + \tau^- \\
{\rm dom}_{\sigma \times \tau}(x) &:\equiv & {\rm dom}_\sigma(\fst \, x) \land {\rm dom}_\tau(\snd \, x) \\
{\rm app}_{\sigma \times \tau}(x,y,z) &:\equiv & \big( \, \exists u^{\sigma^-} \, ( z = \inl \, u \land {\rm app}_\sigma(\fst \, x, \fst \, y, u)) \lor \\ & & \exists v^{\tau^-} \, ( z = \inr \, v \land {\rm app}_\tau(\snd \, x, \snd \, y, v)) \, \big)\\
& & \land {\rm dom}_{\sigma \times \tau}(x) \land {\rm dom}_{\sigma \times \tau}(y) \\ \\
(\sigma + \tau)^+ & :\equiv & \sigma^+ + \tau^+ \\
(\sigma + \tau)^- & :\equiv & \sigma^- \times \tau^- \\
{\rm dom}_{\sigma + \tau}(x) & :\equiv & \exists u^{\sigma^+} (x = \inl \, u \land {\rm dom}_{\sigma}(u)) \lor \exists v^{\tau^+} (x = \inr \, v \land {\rm dom}_{\tau}(v)) \\
{\rm app}_{\sigma + \tau}(x,y,z) & :\equiv & \big( \, \exists u^{\sigma^+}, v^{\tau^+} \, \big( \, ( x = \inl \, u \land y = \inr \, v) \lor (x = \inr \, v \land y = \inl \, u) \big) \lor \\
& & \exists u^{\sigma^+}, v^{\sigma^+} \, ( \, x = \inl \, u \land y = \inl \, v \land {\rm app}_{\sigma^+}(u,v,\fst \, z)) \lor \\
& & \exists u^{\tau^+}, v^{\tau^+} \, ( \, x = \inr \, u \land y = \inr \, v \land {\rm app}_{\tau^+}(u,v,\snd \, z)) \, \big) \land \\
& & {\rm dom}_{\sigma + \tau}(x) \land {\rm dom}_{\sigma + \tau}(y) \\ \\
(\sigma \to \tau)^+ & :\equiv & (\sigma^+ \to \tau^+) \times (\sigma^+ \to \sigma^+ \to \tau^- \to \sigma^-) \\
(\sigma \to \tau)^- & :\equiv & \sigma^+ \times \tau^- \\
{\rm dom}_{\sigma \to \tau}(x) & :\equiv & \forall u^{\sigma^+} ( \, {\rm dom}_\sigma(u) \to {\rm dom}_\tau((\fst \, x)u) \, ) \land \\
& & \forall u^{\sigma^+}, v^{\sigma^+}, w^{\tau^-} \, ( \, {\rm dom}_\sigma(u) \land {\rm dom}_\sigma(v) \land \\ & & {\rm app}_\tau((\fst \, x)(u), (\fst \, x)(v),w) \to  {\rm app}_\sigma(u,v,(\snd \, x)uvw) \, ) \\
{\rm app}_{\sigma \to \tau}(x,y,z) &:\equiv& {\rm dom}_\sigma(\fst \, z) \land {\rm app}_\tau((\fst \, x)(\fst \, z), (\fst \, y)(\fst \, z), \snd \, z) \\
& & \land {\rm dom}_{\sigma \to \tau}(x) \land {\rm dom}_{\sigma \to \tau}(y)
\end{eqnarray*}

Now we the following statements are provable in $\eha$ for each type $\sigma$:
\begin{enumerate}
    \item If ${\rm app}_\sigma(x,y,z)$, then ${\rm dom}_\sigma(x)$ and ${\rm dom}_\sigma(y)$.
    \item If ${\rm dom}_\sigma(x)$, then there is no $z$ such that ${\rm app}_\sigma(x,x,z)$.
    \item There is a functional $s: \sigma^+ \to \sigma^+ \to \sigma^- \to \sigma^-$, such that if ${\rm app}_\sigma(x,y,z)$, then ${\rm app}_\sigma(y,x,sxyz)$.
    \item There is a functional $t: \sigma^+ \to \sigma^+ \to \sigma^+ \to \sigma^- \to (\sigma^- + \sigma^-)$ such that whenever ${\rm app}_\sigma(x,y,u)$ and ${\rm dom}_\sigma(z)$, then either there is some $v$ such that $txyzu = \inl \, v$ and ${\rm app}_\sigma(x,z,v)$ or there is some $w$ such that $txyzu = \inr \, w$ and ${\rm app}_\sigma(y,z,w)$.
\end{enumerate}    

In what follows we assume that we have assigned to each variable $x$ of type $\sigma$ a variable $x^\alpha$ of type $\sigma^+$ in such a way that for any pair of distinct variables $x, y$ of type $\sigma$ the variables $x^\alpha$ and $y^\alpha$ are distinct as well.

\begin{defi} To each term $t$ of $\eha$ of type $\sigma$ with free variables $x_1,\ldots,x_n$ of types $\sigma_1,\ldots,\sigma_n$, respectively, we associate a term $t^\alpha$ of $\eha$ of type $\sigma^+$ with free variables $x_1^\alpha,\ldots,x_n^\alpha$ of types $\sigma^+_1,\ldots,\sigma^+_n$, respectively.
\begin{eqnarray*}
    x^\alpha & :\equiv & x^\alpha \mbox{ if } x \mbox{ is some variable}\\
    (s \, t)^\alpha & :\equiv & (\fst \, s^\alpha) \,t^\alpha
\end{eqnarray*}
The cases of the combinators are quite tedious. We do a few of them which do not involve (too many) nested implications:
\begin{eqnarray*}
    \fst^\alpha & :\equiv & ({\bf fst},\lambda x^{\sigma^+ \times \tau^+}.\lambda y^{\sigma^+ \times \tau^+}.\lambda z^{\sigma^-}.\inl(z)) \\
    \snd^\alpha & :\equiv & ({\bf snd},\lambda x^{\sigma^+ \times \tau^+}.\lambda y^{\sigma^+ \times \tau^+}.\lambda z^{\tau^-}.\inr(z)) \\
    \pair^\alpha & :\equiv & (\lambda x^{\sigma^+}.(\lambda y^{\tau^+}.\pair \, x \, y, \lambda u^{\tau^+}.\lambda v^{\tau^+}.\lambda w^{\sigma^- + \tau^-}.{\bf case} \, (\lambda z^{\sigma^-}.0^{\tau^-}) \, (\lambda z^{\tau^-}.z) \, w), \\ & & \lambda u^{\sigma^+}.\lambda v^{\sigma^+}.\lambda w^{\tau^+ \times (\sigma^- + \tau^-)}.{\bf case} \, (\lambda z^{\sigma^-}.z) \, (\lambda z^{\tau^-}.0^{\sigma^-}) \, (\snd \, w)) \\
    {\bf k}^\alpha & :\equiv & (\lambda x^{\sigma^+}.(\lambda y^{\tau^+}.x, \lambda u^{\tau^+}.\lambda v^{\tau^+}.\lambda w^{\sigma^-}.0^{\tau^-}),\lambda u^{\sigma^+}.\lambda
     v^{\sigma^+}.\lambda w^{\tau^+ \times \sigma^-}.\snd \, w) \\
     & \ldots
\end{eqnarray*}
Here 0 is some closed term of the appropriate type.
\end{defi}

Let us also discuss in some detail how one may interpret the recursor ${\bf R}$. In fact, it will be more convenient to interpret the iterator ${\bf it}_\sigma$ of type $((\sigma \to \sigma) \times \sigma \times 0) \to \sigma$ satisfying the equations:
\begin{eqnarray*}
    {\bf it}_\sigma \, (f,x, 0) & =_\sigma & x, \\
    {\bf it}_\sigma \, (f,x,Sn) & =_\sigma & f({\bf it}_\sigma \, (f,x,n)).
\end{eqnarray*}
Indeed, once such an iterator has been defined, one can define the recursor as follows:
\[ {\bf R}_\sigma = \lambda x^\sigma.\lambda f^{0 \to \sigma \to \sigma}. \lambda n^0.\snd({\bf it}_{0 \times z} \, (\lambda z^{0 \times \sigma}.(S (\fst \, z), f \, (\fst \, z) \, (\snd \, z) ) , (0,x) , n)). \]
The interpretation of the iterator ${\bf it}_\sigma^\alpha = ({\bf pit}_\sigma^\alpha,{\bf nit}_\sigma^\alpha)$ is given by  
\[ {\bf pit}_\sigma^\alpha = \lambda x^{((\sigma^+ \to \sigma^+) \times (\sigma^+ \to \sigma^+ \to \sigma^- \to \sigma^-)) \times (\sigma^+ \times 0)}.{\bf it}_{\sigma^+} (\fst \, (\fst \, x), \fst \, (\snd \, x), \snd \, (\snd \, x)) \]
and ${\bf nit}^\alpha_{\sigma}$ is the following informally described algorithm: suppose we are given $f,x,n$ and $g,y,m$ and evidence showing that ${\bf pit}^\alpha_{\sigma}(f,x,n)$ and ${\bf pit}^\alpha_{\sigma}(g,y,m)$ are apart. Our task is to evidence showing that either $f$ and $g$ are apart, $x$ and $y$ are apart or $n$ and $m$ are apart. We can start by comparing $m$ and $n$: if they are distinct, we provide evidence that they are and we are done. So assume $m = n$. If $m = n = 0$, then we are given evidence that $x$ and $y$ are apart and we are also done. So assume $m = n = Sk$. Then we are given evidence that $f({\bf pit}^\alpha_{\sigma}(f,x,k))$ and $g({\bf pit}^\alpha_{\sigma}(g,y,k))$ are apart. By comparing with $f({\bf pit}^\alpha_{\sigma}(g,y,k))$, we have evidence either that $g({\bf pit}^\alpha_{\sigma}(g,y,k))$ and $f({\bf pit}^\alpha_{\sigma}(g,y,k))$ are apart or $f({\bf pit}^\alpha_{\sigma}(f,x,k))$ and $f({\bf pit}^\alpha_{\sigma}(g,y,k))$ are apart. In the former case, we have evidence that $f$ and $g$ are apart and we are done. In the second case, we can use that $f$ reflects apartness to find evidence that ${\bf pit}^\alpha_{\sigma}(f,x,k)$ and ${\bf pit}^\alpha_{\sigma}(g,y,k)$ are apart, in which case we repeat the procedure until $k = 0$ and we are in one of the previous cases.

The crucial property of this interpretation of the combinators is that, provably in $\eha$, ${\rm dom}({\bf c})$ holds for any combinator $\bf c$ and with respect to the application $s^{(\sigma \to \tau)^+} * t^{\sigma^+} = (\fst \, s) \, t$ they satisfy the usual equations.

\begin{lem}
    If $t$ is a term of $\eha$ of type $\sigma$ and the free variables of the term $t$ are $x_1,\ldots,x_n$ of type $\sigma_1,\ldots,\sigma_n$, respectively, then  
    \[ \eha \vdash \big( \, \bigwedge_i {\rm dom}_{\sigma_i}(x_i^\alpha) \, \big) \to {\rm dom}_{\sigma}(t^\alpha). \] 
\end{lem}
\proof
    By induction on the structure of the term $t$.
\qed

\begin{defi}
    To each formula $\varphi$ of $\eha$ with free variables $x_1,\ldots,x_n$ of types $\sigma_1,\ldots,\sigma_n$, respectively, we associate a formula $\varphi^\alpha$ of $\eha$ with free variables $x_1^\alpha,\ldots,x_n^\alpha$ of types $\sigma^+_n,\ldots,\sigma^+_n$, respectively.
\begin{eqnarray*}
    (s =_\sigma t)^\alpha & :\equiv & \lnot \exists x^{\sigma^-} \, {\rm app}(s^\alpha, t^\alpha, x) \\
    \bot^\alpha & :\equiv & \bot \\
    (\varphi \Box \psi)^\alpha & :\equiv & \varphi^\alpha \Box \psi^\alpha \mbox{ for } \Box \in \{ \lor, \land, \to \} \\
    (\exists x^\sigma \, \varphi)^\alpha & :\equiv & \exists x^{\sigma^+} \, ( \, {\rm dom}_\sigma(x) \land \varphi^\alpha) \\
    (\forall x^\sigma \, \varphi)^\alpha & :\equiv & \forall x^{\sigma^+} \, ( \, {\rm dom}_\sigma(x) \to \varphi^\alpha)
\end{eqnarray*}

\end{defi}

\begin{thm} If $\varphi$ is a formula in the language of $\eha$ with free variables $x_1,\ldots,x_n$ of types $\sigma_1,\ldots,\sigma_n$, respectively, and $\eha + \CE_0 \vdash \varphi$, then \[ \eha \vdash \big( \, \bigwedge_i {\rm dom}_{\sigma_i}(x^\alpha_i) \, \big) \to \varphi^\alpha. \]
\end{thm}
\proof
    By induction on the derivation of $\eha + \CE_0 \vdash \varphi$.
\qed

Very roughly, the reason why the interpretation $\alpha$ witnesses $\CE_0$ is that while elements of type $1 = N \to N$ are interpreted as the elements of type $1$ from the interpreting system (the extra information they are endowed with is vacuous), elements $\Phi$ of type $2 = (N \to N) \to N$ are endowed with the additional information that finds a point $x$ such that $fx \not= gx$ as soon as $\Phi f \not= \Phi g$.

\end{document}